 \documentclass[11pt,twoside]{article}
 \usepackage[dvips]{graphics}
 \usepackage{amsmath,graphicx,amsfonts,amsrefs,amssymb}
 \usepackage{setspace, enumitem, mathtools, float, yfonts}
 \usepackage{graphicx}
  \usepackage[dvips]{graphics}
  \usepackage[section]{placeins}

 \newtheorem{theorem}{Theorem}[section]
 \newtheorem{lemma}[theorem]{Lemma}

 \newtheorem{claim}[theorem]{Claim}

 \def\qed{\hfill\rule{1ex}{1ex}\\}

\begin{document}

 \title{Discussion of a uniqueness result in ``Equilibrium Configurations for a Floating Drop'' }
 \author{Ray Treinen\footnote{Department of Mathematics, Texas State University, 601 University Dr., San Marcos, TX 78666, rt30@txstate.edu}}
 \maketitle

\begin{abstract}
We analyze a uniqueness result presented by Elcrat, Neel, and Siegel \cite{ENS} for unbounded liquid bridges, and show that the proof they presented is incorrect.  We add a choice of three hypothesis to their stated theorem and show that their result holds under this condition.  Then we use Chebyshev spectral methods to build a numerical method to approximate solutions to a related boundary value problems that show one of the three hypothesis holds.
\end{abstract}

\section{Introduction}
\label{intro}

In 2004 Elcrat, Neel, and Siegel published a collection of results on the floating drop problem and the related floating bubble problem \cite{ENS}.  Physically, one can visualize a drop of oil resting on a reservoir of water, and the resulting free boundary problem will not be described in detail here.   This work has been held in high esteem in the field of capillarity, which is evident in the review  Robert Finn wrote in Math Reviews for the paper \cite{Finn}.  We offer a select quote from that review here: 
\begin{quote}
The problem of characterizing the configuration of a drop of liquid floating in equilibrium on the surface of an infinite bath of another liquid appeared initially in the second supplément to the tenth book of Laplace's Mécanique Céleste, in 1806, without detailed treatment. It was later studied by Poisson in his ``Nouvelle Théorie…''(1831), and discussed by Bowditch (1839) in his English translation of Laplace's treatise. These works were remarkable for their time but far from complete, and there appears to be no further mathematical discussion in the literature, prior to the present study. That the problem was ignored so long despite its evident theoretical and also practical interest is perhaps indicative of the technical obstacles that have impeded a full formal description.

Also the present authors have not solved the problem completely, but what they offer is impressive, with a direct hands-on approach. The authors make clever use of new results that appeared in other contexts during the past quarter century.
\end{quote}

Given this perspective, it is unfortunate that the subject of this current  paper is a flaw in one of the proofs presented in that work.  The theorem stated as Theorem~3.2 in \cite{ENS} treats the existence and uniqueness of a boundary value problem for unbounded liquid bridges, and the proof of uniqueness is the topic of discussion here.  We state the theorem here in Theorem~\ref{thm:ENS}.  This flaw is seemingly in the proof alone, and does not seem to be in the stated result.  In what follows we will analyze the presented proof and show how the published approach could potentially be repaired.  Then we offer an alternative approach and we give strong numerical evidence that this alternate approach yields the (still unproven) result found in that paper.  

Before we proceed to the details, a short comment on the impact of this flaw is appropriate.  We note that the strongest results found in \cite{ENS} are not effected by this flaw, but those stronger results do not hold for the general cases of all physical configurations of fluids.  Specifically, they found that under some restrictions on the associated surfaces tensions the floating drop problem is solved for any given drop volume.   This assumption implies all of the component surfaces for a floating drop can be shown to be a graph over a base domain, and under this restriction all of the results in \cite{ENS} still hold.  This assumption is discussed there, and with the goal of avoiding the quite technical description of the floating drop problem, we refer the reader to that work for the detailed criterion.    It is when the problem was generalized to admit all possible physical configurations that the result described below was used.  The results described in that work are Theorem~4.1, Theorem~5.3,  Corollary 5.1 (gaps), and Theorem~4.2 (no gap).

Finally, the layout of this paper is as follows.  In Section~\ref{section:problem} we present the unbounded liquid bridge and the theorem from \cite{ENS}.  In Section~\ref{prelim} we collect preliminary results found in a paper by Vogel \cite{Vogel}.  In Section~\ref{Section:flaw} we analyze the proof given in \cite{ENS}, and in Section~\ref{Section:fix} we prove a theorem on uniqueness with a choice of three hypothesis.  Then in Section~\ref{Section:T} we develop a new Chebyshev spectral method to compute the height of the vertical point of the unbounded liquid bridge as a function of radius.  This function is then used to produce strong numerical evidence that the result in question is true, and the details are found in Section~\ref{Section:rdot}.

\section{Uniqueness of solutions for a boundary value problem involving  unbounded liquid bridges}
\label{section:problem}

\begin{figure}[t]
	\centering
	\scalebox{0.45}{\includegraphics{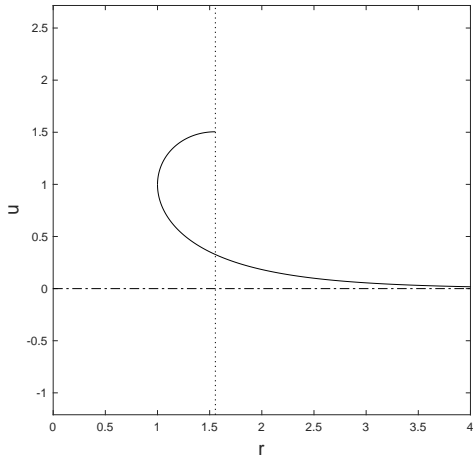}}
	\caption{An example unbounded liquid bridge with vertical point at $\sigma = 1$, where the generating curve for the radially symmetric surface is graphed.  The reference height is indicated at $u =0$ and the vertical line corresponds to the radius where the boundary conditions hold, as discussed in Theorem~\ref{thm:ENS}.}
	\label{fig:unbd_bridge}
\end{figure}

We consider here the family of radially symmetric unbounded liquid bridges.  These bridges are solutions of 
\begin{eqnarray}
	\frac{dr}{d\phi} &=& \frac{-r\cos\phi}{ru + \sin\phi} \label{eqn:drdphi}\\
	\frac{du}{d\phi} &=& \frac{-r\sin\phi}{ru + \sin\phi} \label{eqn:dudphi}
\end{eqnarray}
over the range $\phi \in [0,\pi]$,  with the radial value $r$ and the height of the interface above a fixed reference level is $u$.  The vertical points along the curve are required to be at $(\sigma,T(\sigma))$, where $\phi = \pi/2$.   This condition on the vertical points $(r(\pi/2),u(\pi/2))$ produces solutions that are asymptotic to the $r$-axis.  See Siegel \cite{Siegel} and Vogel \cite{Vogel} for background and more references on the function $T(\sigma)$.  Figure~\ref{fig:unbd_bridge} shows an example liquid bridge of this type.  

The following theorem was stated in \cite{ENS}:
\begin{theorem}\label{thm:ENS}
	For every $\rho_0> 0$ and every $\phi_0 \in [0,\pi)$, there exists a unique symmetric unbounded liquid bridge that meets the radius $r = \rho_0$ so that $\phi = \phi_0$. That is, there exists a unique $\sigma$ satisfying  $r(\phi_0;\sigma) = \rho_0$ so that  solutions of  \eqref{eqn:drdphi}-\eqref{eqn:dudphi} form an unbounded liquid bridge by meeting the requirement that $r(\pi/2;\sigma) = \sigma$ and $u(\pi/2;\sigma) = T(\sigma)$.
\end{theorem} 

We follow the proof presented in \cite{ENS}, noting the flaw in the argument.  Hypotheses will be added later to remedy the flaw, which will be carefully stated in Theorem~\ref{thm:hyp}, and numerical results will be presented that show strong evidence that the general result is true, despite the incorrect proof.

\section{Preliminaries}
\label{prelim}

\begin{figure}[t]
	\centering
	\scalebox{0.45}{\includegraphics{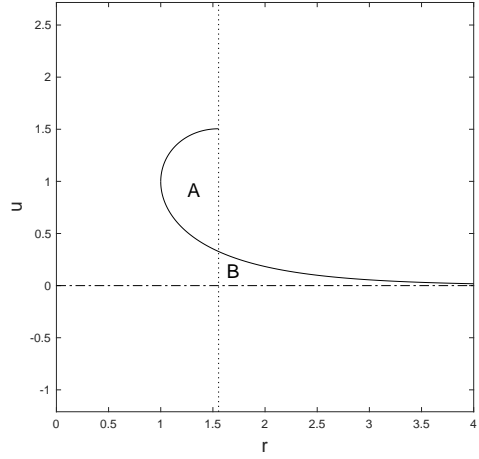}}
	\caption{The regions $A$ and $B$ from Lemma~\ref{lemma:volume}.}
	\label{fig:volume}
\end{figure}

Before we discuss the presented proof of Theorem~\ref{thm:ENS}, we  collect a few results due to Vogel \cite{Vogel}:
\begin{lemma}
	\label{lemma:volume}
Let $\Gamma = (r(\phi),u(\phi))$ be a particular profile curve.  Pick $\phi_0\in[0,\pi)$, and let $r_0 = r(\phi_0)$.  Let $A$ be the solid obtained by rotating the region bounded by $r = r_0$ and $\Gamma$, and let $B$ be the solid obtained by rotating the unbounded region between $\Gamma$ and the $r$-axis from $r=r_0$ to $r=\infty$ around the $u$-axis.  Then
\begin{equation}
|B| - |A| = 2\pi r_0\sin\phi_0.
\end{equation}
\end{lemma}
This lemma is illustrated in Figure~\ref{fig:volume}.
\begin{lemma}
	\label{lemma:equalheight}
Let $\Gamma_1$ and $\Gamma_2$ be two profile curves as above.  If there exists some $\phi_0\in[0,\pi)$ such that  $u_1(\phi_0) = u_2(\phi_0)$, then $\Gamma_1 \equiv \Gamma_2$.
\end{lemma}
\begin{lemma}
No two distinct profile curves can cross twice.
\end{lemma}
The following result is contained in a remark in Vogel's paper, and also is extended  by using Lemma~\ref{lemma:equalheight}
\begin{lemma}\label{monotone}
Given two profile curves with vertical points at radii $\sigma_1$ and $\sigma_2$, if $\sigma_2 > \sigma_1$, then $u(\phi_0;\sigma_2) > u(\phi_0;\sigma_1)$.  Conversely, if $u(\phi_0;\sigma_2) > u(\phi_0;\sigma_1)$, then $\sigma_2 > \sigma_1$.
\end{lemma}
Next we will consider how the system \eqref{eqn:drdphi}-\eqref{eqn:dudphi} behaves as the location of the vertical point moves by differentiating with respect to that parameter $\sigma$, the results of which we denote by $\dot r$ and $\dot u$:
\begin{eqnarray}
\frac{d\dot r}{d\phi} &=& = \cos\phi\frac{\dot u r^2 - \dot r\sin\phi}{(ru + \sin\phi)^2} \label{eqn:dotr}\\
\frac{d\dot u}{d\phi} &=& = \sin\phi\frac{\dot u r^2 - \dot r\sin\phi}{(ru + \sin\phi)^2}.\label{eqn:dotu} 
\end{eqnarray}
The conditions that  $r(\pi/2;\sigma) = \sigma $ and $u(\pi/2,\sigma) = T(\sigma)$ become 
\begin{eqnarray}
\dot r(\pi/2,\sigma) &=& 1, \mbox{ and} \label{eqn:bcdotr}\\
\dot u(\pi/2,\sigma) &=& T^\prime(\sigma). \label{eqn:bcdotu}
\end{eqnarray}
Then Vogel finds
\begin{lemma}
Let $\sigma > 0$ be fixed.  If $\phi_0\in[0,\pi)$ gives $\dot r(\phi_0; \sigma) \geq 0$, then $\dot r(\phi;\sigma) \geq 0$ for $\phi_0 \leq \phi < \pi$.  If  $\dot r(\phi_0; \sigma) > 0$, then $\dot r(\phi;\sigma) > 0$ for $\phi_0 < \phi < \pi$.
\end{lemma}
The second statement follows from the Vogel's same proof.
\begin{lemma}
If  $\dot u(0; \sigma) = 0$, then $\dot r(0;\sigma) < 0$.
\end{lemma}
\begin{lemma}
If $\phi \in(0,\pi)$, then $\dot u(\phi;\sigma) > 0$.
\end{lemma}

\section{The proof presented}
\label{Section:flaw}

The proof is by contradiction.  Let $\Gamma_1$ and $\Gamma_2$ be two profile curves as above such that $r_1(\phi_0) = r_2(\phi_0)$ for some $\phi_0\in[0,\pi/2)$.  We will assume that these curves exist and are distinct.  Upon possibly relabeling, we have $u_1(\phi_0) > u_2(\phi_0)$.

Suppose that $\phi_0$ is the largest such value of $\phi$ such that $r_1(\phi_0) = r_2(\phi_0)$ and denote that radius by $\rho_0$.  These are the leftmost such points, and $\rho_0$ is the smallest such radius.  Let $A_1$, $B_1$, $A_2$, and $B_2$ be the corresponding volumes from Lemma~\ref{lemma:volume} for these profile curves.

\begin{figure}[h]
	\centering
	\scalebox{0.45}{\includegraphics{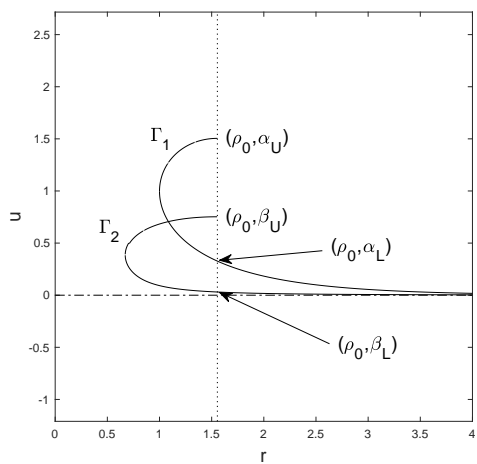}}
	\caption{The two (potentially) distinct curves $\Gamma_1$ and $\Gamma_2$.}
	\label{fig:alphabeta}
\end{figure}

Consider the intersection points of these curves with $r = \rho_0$.  Define $\alpha_U = u_1(\phi_0)$ and $\beta_U = u_2(\phi_0)$ to be the upper intersection points, and define $\alpha_L$ and $\beta_L$ to be the lower intersection points.   See Figure~\ref{fig:alphabeta}.  We will have some cases to consider, depending on these values.

{\bf Step 1} (The case that $\alpha_L \leq \beta_L$)
This means that $\alpha_U > \beta_U > \beta_L \geq \alpha_L$.  Since $\Gamma_1$ and $\Gamma_2$ cannot cross more than once, in this case they cannot cross at all.  Thus $\Gamma_1$ lies entirely outside the region bounded by $\Gamma_2$ and the line $r = \rho_0$.  Thus $\sigma_1 < \sigma_2$ and it follows that $T(\sigma_1) < T(\sigma_2)$.  This then implies that there is some $\phi_1\in(\phi_0,\pi/2)$ such that $u_1(\phi_1) = u_2(\phi_1)$.  This contradicts Lemma~\ref{lemma:equalheight} and eliminates the case that $\alpha_L \leq \beta_L$..

{\bf Step 2} (The case that $\alpha_L > \beta_L$) 
As above, if $T(\sigma_1) \leq T(\sigma_2)$, we contradict  Lemma~\ref{lemma:equalheight}.  Thus $T(\sigma_1) > T(\sigma_2)$ and $\sigma_1 > \sigma_2$.  The goal here is to contradict Lemma~\ref{lemma:volume}.  To this end, we consider two more steps, the first of which is to show that $|B_1| > |B_2|$ and the second is to show that $|A_1| < |A_2|$.

{\bf Step 3 }
If $\alpha _L \leq \beta_U$, then $\Gamma_1$ and $\Gamma_2$ must cross somewhere above  $\alpha_L$.  Then, as they cannot cross twice, they cannot cross below $\alpha_L$.  Then, as $\alpha_L  > \beta_L$, this implies that the lower portion of $\Gamma_1$ lies completely above the lower portion of $\Gamma_2$ when $ r >  \rho_0$.  Thus $\alpha _L \leq  \beta_U$ implies $|B_1| > |B_2|$.

If $\alpha _L \geq \beta_U$, then the part of $\Gamma_1$ with $r < \rho_0$ lies completely above the corresponding part of $\Gamma_2$.  Thus, at $u = \beta_U$, $\Gamma_1$ has an inclination angle greater than $\pi/2$, while $\Gamma_2$ has an inclination angle less than $\pi/2$.  Also, at $r = \rho_0$, the lower part of $\Gamma_1$ is above the lower part of $\Gamma_2$.  the goal is to show that $\Gamma_1$ remains above $\Gamma_2$ for all $r > \rho_0$.  Assume the opposite: the two curves cross somewhere below $u = \beta_U$.  At this crossing point the inclination angle of $\Gamma_2$ would be less than that of $\Gamma_2$.  At $u = \beta_U$ the inequality was reversed.  This implies that there is a value of $u$ between $\beta_U$ and the value of $u$ at the crossing point where both curves have the same inclination angle.  This contradicts Siegel's uniqueness theorem \cite{Siegel}, which would imply that $\Gamma_1 \equiv \Gamma_2$ when that happens.  Thus the curves cannot cross, and the lower part of  $\Gamma_1$ lies above the lower part of $\Gamma_2$ for all $r > \rho_0$.  Thus  $\alpha _L \geq \beta_U$ implies $|B_1| > |B_2|$, and we  have established that  $|B_1| > |B_2|$ holds for all possibilities in the case that $\alpha_L > \beta_L$.

{\bf Step 4 } (Incorrect)
Let $\Gamma^\prime_1$ be the rigid translation of $\Gamma_1$ downward by $\alpha_U - \beta_U$ so that $\Gamma^\prime_1$ and $\Gamma_2$ are tangent at their topmost points.  We are concerned with the regions bounded by $\Gamma^\prime_1$ and $r = \rho_0$ and by $\Gamma_2$ and $r = \rho_0$.  We seek to show that $\Gamma^\prime_1$ is contained in the region bounded by $\Gamma_2$.

Starting at $(\rho_0,\beta_U)$, both profile curves are functions $u(r)$, so that we may consider the second derivative:
\begin{eqnarray}
\frac{d^2u}{dr^2} &=& \frac{d\,\,\,}{dr} \tan\phi \nonumber\\
&=& \sec^2\phi \frac{d\phi}{dr} \nonumber\\
&=& \frac{1}{\cos^2\phi}\frac{ru + \sin\phi}{-r\cos\phi}\nonumber \\
&<& 0,
\end{eqnarray}
which corrects a minor (inconsequential) error in \cite{ENS}.

At $\phi = \phi_0$, we have $r = \rho_0$, but $u_1 > u_2$ for $\Gamma^\prime_1$ and $\Gamma_2$ where we draw values of $u$ from the underlying equations \eqref{eqn:drdphi}-\eqref{eqn:dudphi}.  Thus, as  $\frac{d^2u}{dr^2} < 0$, on $\Gamma^\prime_1$ we have a smaller value for $\frac{d^2u}{dr^2}$, and at least on some small interval to the left of  $r = \rho_0$, $\Gamma^\prime_1$ lies below $\Gamma_2$, and has a greater inclination angle.  Thus we know that $\Gamma^\prime_1$ starts ``inside'' $\Gamma_2$.  Is it possible it can ever leave? 

Suppose that $\Gamma^\prime_1$ leaves on the upper portion of $\Gamma_2$.  Then at this point of leaving, $\Gamma^\prime_1$ has a smaller inclination angle than that of $\Gamma_2$,  This implies that there exists some $r^*$ between the crossing point and $r=\rho_0$ where both $\Gamma^\prime_1$ and $\Gamma_2$ have the same inclination angle.  That means, of course, that both $\Gamma_1$ and $\Gamma_2$ have the same inclination angle at $r^* < \rho_0$.  This contradicts the definition of $\rho_0$ to be the smallest such radius with the same inclination angle.  We conclude that $\Gamma^\prime_1$ cannot leave across the upper portion of $\Gamma_2$.

The next part is to eliminate the case where $\Gamma^\prime_1$ leaves across the lower portion of $\Gamma_2$ while $r \leq\rho_0$.  If $\Gamma^\prime_1$ did escape in this region of $\Gamma_2$, it would not be able to return.  The contradiction in that situation is similar to the argument in the previous paragraph.  If $\Gamma^\prime_1$ escaped on the lower portion of $\Gamma_2$, then it would have a smaller inclination angle than that of $\Gamma_2$ there.  Thus, somewhere in between there would be a radius $r^* < \rho_0$ where $\Gamma_1$ and $\Gamma_2$ have the same inclination angle, which is a contradiction.

This means that a necessary condition for $\Gamma^\prime_1$ to escape $\Gamma_2$ is that 
$$
\alpha_U - \alpha_L > \beta_U - \beta_L.
$$
We observe that 
\begin{eqnarray}
\alpha_U - \alpha_L &=& (\alpha_U - T(\sigma_1)) + (T(\sigma_1) - \alpha_L), \nonumber\\
\beta_U - \beta_L &=& (\beta_U - T(\sigma_2)) + (T(\sigma_2) - \beta_L),\nonumber
\end{eqnarray}
and we compute
\begin{eqnarray}
\alpha_U - T(\sigma_1) &=& \int^{\pi/2}_{\phi_0} \frac{r_1\sin\phi}{r_1u_1 + \sin\phi}\, d\phi, \nonumber\\
\beta_U - T(\sigma_2) &=& \int^{\pi/2}_{\phi_0} \frac{r_2\sin\phi}{r_2u_2 + \sin\phi}\, d\phi. \nonumber
\end{eqnarray}
Then we compare the integrands at $\phi_0$, and our assumption that $\alpha_U > \beta_U$ means that  $u_1 > u_2$ there, and we find
\begin{equation}
\frac{\rho_0\sin\phi_0}{\rho_0u_1 + \sin\phi_0} < \frac{\rho_0\sin\phi_0}{\rho_0u_2 + \sin\phi_0}.
\end{equation}
\begin{claim}
We have 
\begin{equation}
\frac{r_1\sin\phi}{r_1u_1 + \sin\phi} < \frac{r_2\sin\phi}{r_2u_2 + \sin\phi}
\end{equation}
on the whole interval of integration $(\phi_0,\pi/2)$.
\end{claim}
We will see that the proof provided for this claim is incorrect.  We proceed with the argument provided.

The approach is by contradiction.  If the claim does not hold over the entire interval, then there would be some $\phi$ such that 
\begin{equation}
	\frac{r_1\sin\phi}{r_1u_1 + \sin\phi} = \frac{r_2\sin\phi}{r_2u_2 + \sin\phi},
\end{equation}
giving
\begin{eqnarray}
r_1r_2u_2 + r_1\sin\phi &=& r_1r_2u_1 + r_2\sin\phi, \\
r_1r_2(u_2 - u_1) &=& (r_2 - r_1)\sin\phi.
\end{eqnarray}
The left side is negative throughout the interval, since we are on the upper portions of both curves, and $\Gamma^\prime_1$ cannot leave on the upper part of $\Gamma_2$.  That we are on the upper portions of both of the curves, and that $\Gamma^\prime_1$ has not left the upper part of $\Gamma_2$ implies that $r_2 < r_1$ in this region, and the right hand side is also negative.  The authors incorrectly assert that the right hand side is positive, and draw an incorrect contradiction.

Thus this claim is not proved.

{\bf Step 4 } (Conclusions)

The remainder of the proof presented seems to be without error, however, it depends on the claim holding.  The authors' incorrect argument leads to the conclusion that  
\begin{equation}
|A_1| < |A_2|
\end{equation}
and thus
\begin{equation}
|B_1| - |A_1| > |B_2| - |A_2|,
\end{equation}
which contradicts Lemma~\ref{lemma:volume}, which says
$$
|B_1| - |A_1| = 2\pi\rho_0\sin\psi_0 = |B_2| - |A_2|.
$$
This would have completed the case that $\alpha_L > \beta_L$.

If any approach yields $|A_1| < |A_2|$ in some form, then the approach discussed here  (developed by Elcrat, Neel, and Siegel) yields the uniqueness result under discussion.

\section{Adapting the  approach}
\label{Section:fix}

We then consider what happens to the region $A$ when the initial height is changed.  We will need the formulation of the ODE given by
\begin{equation}
(r\sin\phi)_r = -ru.
\end{equation}
We will also need to precisely describe how the system \eqref{eqn:drdphi}-\eqref{eqn:dudphi} behaves as $u_0$ varies by differentiating with respect to $u_0$, the result of which we denote by $\acute r$ and $\acute u$.  We find
\begin{eqnarray}
	\frac{d\acute r}{d\phi} &=& = \cos\phi\frac{\acute u r^2 - \acute r\sin\phi}{(ru + \sin\phi)^2} \label{eqn:acuter}\\
	\frac{d\acute u}{d\phi} &=& = \sin\phi\frac{\acute u r^2 - \acute r\sin\phi}{(ru + \sin\phi)^2}, \label{eqn:acuteu} 
\end{eqnarray}
and the conditions that  $r(0;u_0) = \rho_0 $ and $u(0,u_0) = u_0$ become 
\begin{eqnarray}
	\acute r(0,u_0) &=& 0, \mbox{ and} \label{eqn:bcacuter}\\
	\acute u(0,u_0) &=& 1. \label{eqn:bcacuteu}
\end{eqnarray}

We compute the volume $V$ of this region directly, using the method of washers to find
\begin{eqnarray}
V(\phi;u_0) &=& \pi\rho^2_0(u_0 - u(\phi)) - \pi\int^{u_0}_{u(\phi)} r^2\, du \nonumber\\
&=& \pi\rho^2_0(u_0 - u(\phi)) - \pi\int^{\phi_0}_\phi r^2\frac{du}{d\phi}\, d\phi \nonumber\\
&=& \pi \rho^2_0(u_0 - u(\phi)) - \pi\left( \rho^2_0 u_0 - r^2 u - 2\int^{\phi_0}_\phi r \frac{dr}{d\phi} u \, d\phi    \right) \nonumber\\
&=& \pi(r^2 - \rho_0^2)u(\phi) + 2\pi\left(  \int^{\rho_0}_\sigma ru\, dr + \int^\sigma_r ru\, dr   \right) \nonumber\\
&=& \pi(r^2 - \rho_0^2)u(\phi) - 2\pi\left(  \int^{\rho_0}_\sigma (r\sin\phi)_r\, dr + \int^\sigma_r (r\sin\phi)_r\, dr   \right) \nonumber\\
&=& \pi(r^2 - \rho_0^2)u(\phi) + 2\pi\left( r\sin\phi - \rho_0\sin\phi_0  \right). 
\end{eqnarray}

Then, as there is some $\phi^- = \phi^-(u_0)> \pi/2$ where $r(\phi^-) = \rho_0$, we find
\begin{equation}
|A| = V(\phi^-;u_0) = 2\pi\rho_0\left( \sin\phi^- - \sin\phi_0  \right).
\end{equation}

We next consider how this changes as we change the initial height.  We denote the derivative with respect to the parameter $u_0$ by an accent so that
\begin{eqnarray}
\acute V &=& 2\pi r \acute r u + \pi\left( r^2 - \rho_0^2   \right)\acute u + 2\pi \acute r \sin\phi \nonumber\\
&=& 2\pi\acute r(ru + \sin\phi) + \pi \left( r^2 - \rho^2_0  \right)\acute u \nonumber\\
&=& 2\pi \acute r \Delta + \pi \left( r^2 - \rho^2_0  \right)\acute u, \nonumber
\end{eqnarray}
where $\Delta :=  ru + \sin\phi$ is positive.  Then, to see how $A$ changes with $u_0$, we evaluate this expression at $\phi^-$, where $r = \rho_0$, to find
\begin{equation}
\acute V = 2\pi\acute r \Delta.
\end{equation}
Thus the sign of $\acute r(\phi^-)$ determines the sign of $\acute V$, and the uniqueness result will follow if $\acute r(\phi^-) < 0$.

It is  sometimes convenient to work with  the parameter $\sigma$ instead of $u_0$.  By Lemma~\ref{monotone} and the analytic dependence of our solutions on the data, we have 
$$
\frac{du_0}{d\sigma} > 0.
$$
Then, using the dot notation for derivatives with respect to the parameter $\sigma$,
\begin{equation}
\dot V(\phi;u_0(\sigma)) = 2\pi \dot r\Delta + \left( r^2 - \rho^2_0    \right)\dot u
\end{equation}
and
\begin{equation}
\dot V(\phi; u_0(\sigma)) = \acute V(\phi; u_0)\frac{du_0}{d\sigma}.
\end{equation}
Thus, at $\phi = \phi^-$, we have
\begin{equation}
\acute V \frac{du_0}{d\sigma} = \dot V = 2\pi\dot r \Delta,
\end{equation}
and the uniqueness result will follow if $\dot r(\phi^-) < 0$.

We will also take an independent approach to proving the desired uniqueness result that does not depend on the approach used by Elcrat, Neel, and Siegel.  Denote $(r,u)$ to be solutions of \eqref{eqn:drdphi}-\eqref{eqn:dudphi}  that meets the radius $r = \rho_0$ so that $\phi = \phi_0$.  First, we present a short argument for the existence of solutions.   We know from Vogel \cite{Vogel} that the vertical points $(\sigma,T(\sigma))$ parameterized by $\sigma$ describes all such surfaces, and that this function $T(\sigma)$ is differentiable.  Vogel shows that 
$$
\lim_{\sigma \rightarrow 0} \frac{u(0;\sigma)}{r(0;\sigma)} = 0,
$$
which, with the estimate
$$
r(0;\sigma) = \mathcal{O}\left(   \frac{1}{\sqrt{\log(1/\sigma)}}  \right)
$$
means $u(0;\sigma) \rightarrow 0$ as $\sigma \rightarrow 0$.  Vogel also shows
$$
\sqrt{\frac{\sigma}{T(\sigma)} + \sigma^2} \leq r(0;\sigma) \leq\sqrt{\frac{2\sigma}{T(\sigma)} + \sigma^2},
$$
which implies that $r(0;\sigma) < \infty$ for $\sigma < \infty$.  We also have a bound on $T$ by $\sqrt{2}$.

Then we start with a value of $\sigma$ small enough that $r(0;\sigma) < \rho_0$.  We then smoothly increase $\sigma$ so that $r(0;\sigma_a) = \rho_0$ for some $\sigma_a\in(0,\rho_0)$ and then we continue up to $r(0;\rho_0)$, where $\sigma = \rho_0$ and the entire upper portion of the curve has passed through the radius $\rho_0$.  Thus, by the intermediate value theorem, there exists a $\sigma_0 \in[\sigma_a,\rho_0]$ such that $r(\phi_0;\sigma_0) = \rho_0$.  This establishes the existence of a solution, which we present here for completeness.

Next, we define $F(\phi;\sigma) = r(\phi;\sigma) - \rho_0.$  We have shown that $F(\phi_0;\sigma) = r(\phi_0;\sigma) - \rho_0$ has at least one zero.  Let $\dot r$ be a component of the solution to \eqref{eqn:dotr}-\eqref{eqn:dotu} subject to \eqref{eqn:bcdotr}-\eqref{eqn:bcdotu}.   If $\dot r(\phi_0;\sigma) > 0$ for all $\sigma>0$, then uniqueness would follow.

We have shown
\begin{theorem}
\label{thm:hyp}
Let $\rho_0> 0$ and $\phi_0 \in [0,\pi)$ be fixed and arbitrary.  If we denote $(r,u)$ to be solutions of \eqref{eqn:drdphi}-\eqref{eqn:dudphi}  that meets the radius $r = \rho_0$ so that $\phi = \phi_0$, then these solutions exist.
Denote by $\phi^-\in(\pi/2,\pi)$ the angle where $r(\phi^-) = \rho_0.$
Let one of the following three conditions hold:
\begin{enumerate}
	\item If $\dot r$ is the component of a solution to \eqref{eqn:dotr}-\eqref{eqn:dotu} subject to \eqref{eqn:bcdotr}-\eqref{eqn:bcdotu}, then $\dot r(\phi^-)< 0$.
	\item If $\acute r$ is the component of a solution to \eqref{eqn:acuter}-\eqref{eqn:acuteu} subject to \eqref{eqn:bcacuter}-\eqref{eqn:bcacuteu}, then $\acute r(\phi^-)< 0$.
	\item If $\dot r$ is the component of a solution to \eqref{eqn:dotr}-\eqref{eqn:dotu} subject to \eqref{eqn:bcdotr}-\eqref{eqn:bcdotu}, then $\dot r(\phi_0)> 0$.
\end{enumerate}
Then there exists a unique symmetric unbounded liquid bridge that meets the radius $r = \rho_0$ so that $\phi = \phi_0$.  That is, there exists a unique $\sigma$ satisfying   $r(\phi_0;\sigma) = \rho_0$ so that  solutions of  \eqref{eqn:drdphi}-\eqref{eqn:dudphi} form an unbounded liquid bridge by meeting the requirement that $r(\pi/2;\sigma) = \sigma$ and $u(\pi/2;\sigma) = T(\sigma)$.
\end{theorem} 

This author has not found a way to rigorously establish any of those three conditions.  In the course of the rest of this paper we will give numerical evidence that the third criteria holds.

\section{Computation of $T(\sigma)$ using Chebyshev Spectral Methods}
\label{Section:T}

Given a radius $\sigma > 0$, we will find the height $T(\sigma)$ of the vertical point on the unbounded liquid bridge there.  In order to find this height, we will need to solve \eqref{eqn:drdphi}-\eqref{eqn:dudphi} for $\phi\in[\pi/2, \pi)$ so that the solution has the required height  decay at infinity.  We will adapt a recently developed Chebyshev spectral method to achieve this.
We summarize the basic approach found in \cite{Treinen}, and refer to that work for technical details.  The goal here is to include only enough details of the work in \cite{Treinen} to explain the basic ideas and how they are modified to the present problem.

The equations \eqref{eqn:drdphi}-\eqref{eqn:dudphi} can be written as a system of three nonlinear ordinary differential  equations, parametrized by the arclength $s$:
\begin{eqnarray}
	\frac{dr}{ds} &=& \cos\psi, \label{drds}\\
	\frac{du}{ds} &=& \sin\psi, \label{duds}\\
	\frac{d\psi}{ds} &=& \kappa u - \frac{\sin\psi}{r}, \label{dpsids}
\end{eqnarray}
where we still have $r$ as the radius and $u$ as the height of the interface, and we introduce the inclination angle $\psi$, which merely satisfies $\psi = \phi - \pi$.  
We will specify boundary conditions by the requirement that at some arclength $\ell>0$ the radius $r(\ell)$ meets a prescribed value $b>0$, and the inclination angle $\psi(\ell)$ meets a prescribed value $\psi_b \in [-\pi, \pi]$.  However, this value of the arclength $\ell$ is unknown, so we rescale the problem.  We define $\tau = s/\ell$, or $s = \ell\tau$.  Then we define
\begin{eqnarray}
	R(\tau) &:=& r(\ell\tau) = r(s),  \nonumber\\
	U(\tau) &:=& u(\ell\tau) = u(s), \nonumber\\
	\Psi(\tau) &:=& \psi(\ell\tau) = \psi(s). \nonumber 
\end{eqnarray}
Then, using the chain rule and multiplying each equation by $\ell$, \eqref{drds}-\eqref{dpsids} become
\begin{eqnarray}
	R^\prime(\tau) - \ell\cos\Psi(\tau) &=& 0, \label{eqn:R}\\
	U^\prime(\tau) - \ell\sin\Psi(\tau) &=& 0, \label{eqn:U}\\
	\Psi^\prime(\tau) + \frac{\ell\sin\Psi(\tau)}{R(\tau)} - \kappa\ell U(\tau) &=& 0. \label{eqn:Psi}
\end{eqnarray}
If we define the column vector $\mathbf{v} = [R\,\, U\,\, \Psi\,  \ell]^T$, we can use \eqref{eqn:R}-\eqref{eqn:Psi} to define the nonlinear operator in the vector equation
\begin{equation}
	\tilde N(\mathbf{v}) = \mathbf{0}. \label{eqn:tN=0}
\end{equation}
We then use the boundary conditions
\begin{eqnarray}
	R(1) - b &=& 0, \label{eqn:bcRb} \\
	R(-1) - \sigma &=& 0, \label{eqn:bcRmb} \\
	\Psi(1) - \psi_b &=& 0, \label{eqn:bcPsib} \\
	\Psi(-1) + \pi/2 &=& 0, \label{eqn:bcmPsib}
\end{eqnarray}
so that we have some form of a two-point boundary value problem.  Here we are specifying the vertical point at the prescribed radius $\sigma$, and we postpone our discussion of the conditions at $r = b$ for later.   We append \eqref{eqn:tN=0} with these  boundary conditions to form the system 
\begin{equation}
	{N}(\mathbf{v}) = \mathbf{0}.
\end{equation}

We will approach this nonlinear problem with a Newton method, and we will need to use the Fr\'{e}chet derivative 
$$
F(\mathbf{v}) = \frac{d N}{d\mathbf{v}}(\mathbf{v}).
$$
Since $\mathbf{v}$ has several components, and some of the corresponding components of $N(\mathbf{v})$  involve applying derivatives with respect to $\tau$, we introduce the differential operator 
$$
D = \frac{d\, }{d\tau}, 
$$
which is applied in a block fashion to $\mathbf{v}$ so that $R^\prime(\tau) = [D\,\, 0\,\, 0\,\, 0]\mathbf{v}$, for example.  We will also have need to use an operator version of function evaluation.  We denote $D^0_\tau$ to be this operator, so that $D^0_1R = R(1)$.  With this in hand, we compute
\begin{equation} \label{eqn:N=0}
	F(\mathbf{v} ) = 
	\begin{bmatrix}
		D & 0 & \ell\sin\Psi & -\cos\Psi \\
		0 & D & -\ell\cos\Psi & -\sin\Psi \\
		\frac{-\ell\sin\Psi }{R^2} & -\kappa\ell & D + \frac{\ell\cos\Psi }{R} & \frac{\sin\Psi}{R} - \kappa U \\
		D^0_{-1} & 0 & 0 & 0 \\
		D^0_{1} & 0 & 0 & 0 \\
		0 & 0 & D^0_{-1} &  0 \\
		0 & 0 & D^0_{1} &  0 
	\end{bmatrix}
	\mathbf{v}.
\end{equation}
We will have need to solve linear systems based on the definition $F(\mathbf{v}) := L\mathbf{v}$.

The discretization of the differential operators is in the block form,  as is described in Driscoll and Hale \cite{DriscollHale2016} and then Aurentz and Trefethen \cite{AurentzTrefethen2017} and for related capillary problems in \cite{Treinen}.  Chebyshev differentiation matrices can be realized as the linear transformation between the data points corresponding to the interpolating polynomials for a function $f$ and its derivative $f^\prime$, where the data is sampled at Chebyshev grid points $x_j = \cos(\theta_j)\in[-1,1]$ where the angles $\theta_j$ are equally spaced angles over $[0,\pi]$.  The basic building blocks of  the nonlinear equation are based on $D^0$ and $D$, which we implement using the Chebfun \cite{Chebfun} commands
\begin{verbatim}
	D0 = diffmat([n-1 n],0,X);
	D1 = diffmat([n-1 n],1,X);
\end{verbatim}
where $X = [-1,1]$,  $n$ is the number of Chebyshev points we are using, and the input 0 or 1 indicates the number of derivatives.  Since $\texttt{D0}$  is rectangular, it becomes a $(n-1)\times n$ identity matrix interpreted as a dense ``spectral down-sampling'' matrix implemented as interpolating on an $n$-point grid followed by sampling on an $(n-1)$-point grid.    We sparsely build $N$ and $L$ using these components.  The basic loop is
\begin{verbatim}
while res_newton > tol_newton
      dv = L(v)\N(v);  
      v = v - dv;
      res_newton = norm(dv,'fro')/norm(v,'fro');
end
\end{verbatim}
where we leave a large number of technical details of the adaptive algorithm for further reading in \cite{Treinen}.  We do modify the initial guess for the starting point of Newton's method from that work.  For $\tau \in X$ at values of $x_j$, we define our initial guesses as
\begin{eqnarray}
    R_0(\tau) &=& (1 + \tau)b/2 + (1 - \tau)\sigma/2 \nonumber\\
U_0(\tau) &=&   \exp(-R_0(\tau) + \sigma) \nonumber\\
\Psi_0(\tau) &=&  \tan^{-1}(-\exp(-R_0(\tau) + \sigma)) \nonumber\\
\ell_0 &=&  b - a. \nonumber
\end{eqnarray}

\begin{figure}[t]
	\centering
	\scalebox{0.45}{\includegraphics{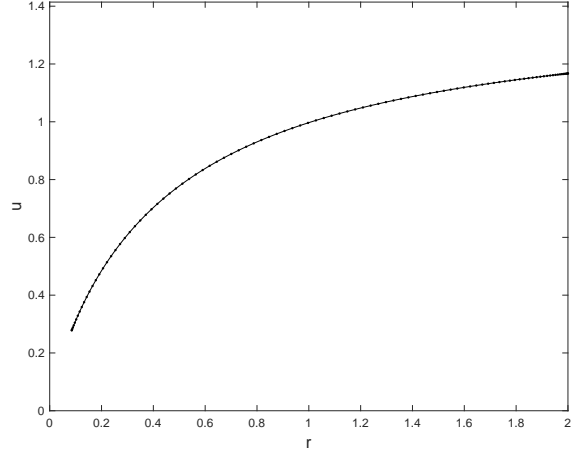}}
	\caption{The function $T(\sigma)$.  The dots along the curve indicate the Chebyshev points and the curve between the dote is the interpolating polynomial, in contrast with the more commonly used linear interpolation between data points.}
	\label{fig:Tcheb}
\end{figure}

\begin{figure}[t]
	\centering
	\scalebox{0.45}{\includegraphics{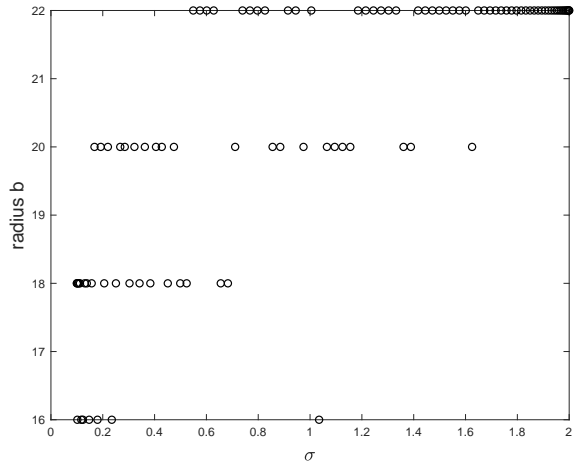}}
	\caption{The ending radius $b$ as a function of the radius $\sigma$ at which $T(\sigma)$ was computed within the requested tolerance. }
	\label{fig:bvec}
\end{figure}

Then, as $T(\sigma)$ is differentiable, it is clearly continuous, so we solve the above boundary value problem with  $b = \max(14,\sigma + 4)$ and $\psi_b = 0$, and then we denote the height of $u$ at $\sigma$ by $T_0$.  Then we increase the value of $b$ by two and we resolve the above boundary value problem with this new $b$ and compare $T_0$ with the new height $T$ of $u$ at $\sigma$.  If the difference $|T - T_0|$  of these heights is within a prescribed absolute error tolerance \texttt{tol\_abs = 1e-11;}, then we have found a reasonable approximation of  $T(\sigma)$, and otherwise we increase $b$ in an iterative fashion until we meet this requested tolerance.  Of course, we reuse the converged data at each step of $b$ to guide the initial guess at the next step for $b$, where we scale the values of $R > 1$ to match at the new value of $b$.  We also used the data in Figure~\ref{fig:bvec} to tune our initial pick of $b$ as 14.

We have released the matlab code for finding $T(\sigma)$ under an open source license and hosted it on a software repository found at \\
\texttt{https://github.com/raytreinen/Unbounded-Liquid-Bridges.git}

\begin{figure}[h]
	\centering
	\scalebox{0.45}{\includegraphics{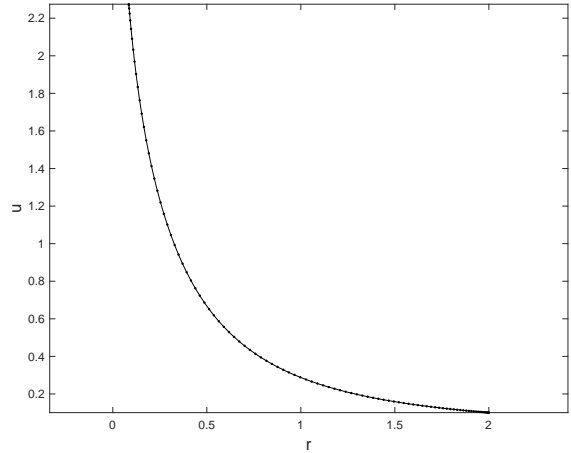}}
	\caption{The function $T^\prime(\sigma)$.}
	\label{fig:Tprime}
\end{figure}

In light of the initial condition \eqref{eqn:bcdotu}, we need a specialized approach in our approximation of $T^\prime(\sigma)$.  We again use a Chebyshev differentiation matrix for good accuracy.  We take 100 Chebyshev points from a positive prescribed $\sigma_{\min}$ to a larger prescribed $\sigma_{\max}$ and we then find $T(\sigma)$ at these Chebyshev points.  We plot the results  in Figure~\ref{fig:Tcheb} for  $\sigma_{\min} = 0.085$ and $\sigma_{\max} = 2$.  For smaller values of $\sigma_{\min}$ we have found that the adaptive algorithm adapted here from \cite{Treinen} uses quite a few Chebyshev points and slows down dramatically.  A multi-scale algorithm would improve this performance, and we leave this for a further work.

Then we build the Chebyshev differentiation matrix for that range of $\sigma$ by using the Chebfun command with \texttt{bigN}$=100$
\begin{verbatim}
diffmat(bigN, 1, [sig_min;sig_max]);
\end{verbatim}
and we plot the result of computing $T^\prime(\sigma)$ in Figure~\ref{fig:Tprime}.  The code for this step is also provided on the repository mentioned above.

\section{Numerical study of $\dot r$}
\label{Section:rdot}

\begin{figure}[b!]
	\centering
	\scalebox{0.45}{\includegraphics{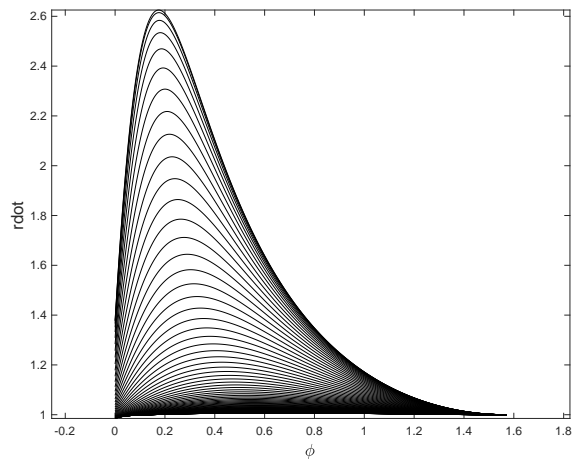}}
	\caption{The solutions $\dot r(\phi)$ foliate a region of the $\phi \dot r-$plane and do not enter the lower half-plane.}
	\label{fig:rdotfoliation}
\end{figure}

\begin{figure}[h]
	\centering
	\scalebox{0.45}{\includegraphics{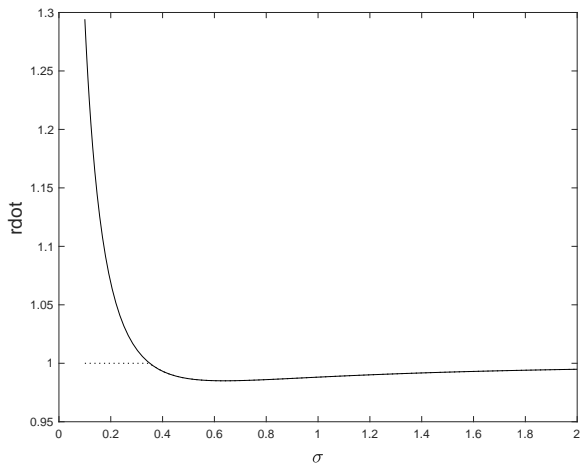}}
	\caption{The endpoint $\dot r(0,\sigma)$ is graphed as a function of $\sigma$.  This is the lowest part of the curve $\dot r(\phi,\sigma)$ in many cases, but for smaller values of $\sigma$ the lowest part of the curve occurs at $\phi = \pi/2$, and this is graphed with a dotted line when relevant. }
	\label{fig:minrdot}
\end{figure}

Now we have the components necessary to construct numerical approximations of the solutions of \eqref{eqn:dotr}-\eqref{eqn:dotu} with the boundary conditions \eqref{eqn:bcdotr}-\eqref{eqn:bcdotu} where we need to append that system with the original \eqref{eqn:drdphi}-\eqref{eqn:dudphi} where $r(\pi/2,\sigma) = \sigma$ and $u(\pi/2,\sigma) = T(\sigma)$.  Then we use Matlab's ode45 to solve this system from $\phi = \pi/2$ to $\phi= 0$ to get the ``top'' portion of the solution.  Here we  ask for 11 digits of accuracy in both the absolute and relative error.  We then use all 100 values of $\sigma$ considered in the last section and sweep out a region of the $\phi \dot r-$plane in Figure~\ref{fig:rdotfoliation}.  This numerical experiment shows that the entirety of the foliation lies completely in the upper half-plane, and thus $\dot r$ is never negative.    We also visualize this phenomenon by graphing the endpoint $\dot r(0,\sigma)$  as a function of $\sigma$.  For many cases, this endpoint is the  lowest part of the curve $\dot r(\phi,\sigma)$, however for smaller values of $\sigma$ the lowest part of the curve occurs at $\phi = \pi/2$, and this is also included in Figure~\ref{fig:minrdot}.   This second figure further illustrates the conclusion of this numerical experiment.  We also included the  code for these experiments on the repository mentioned above.

\begin{figure}[t]
	\centering
	\scalebox{0.45}{\includegraphics{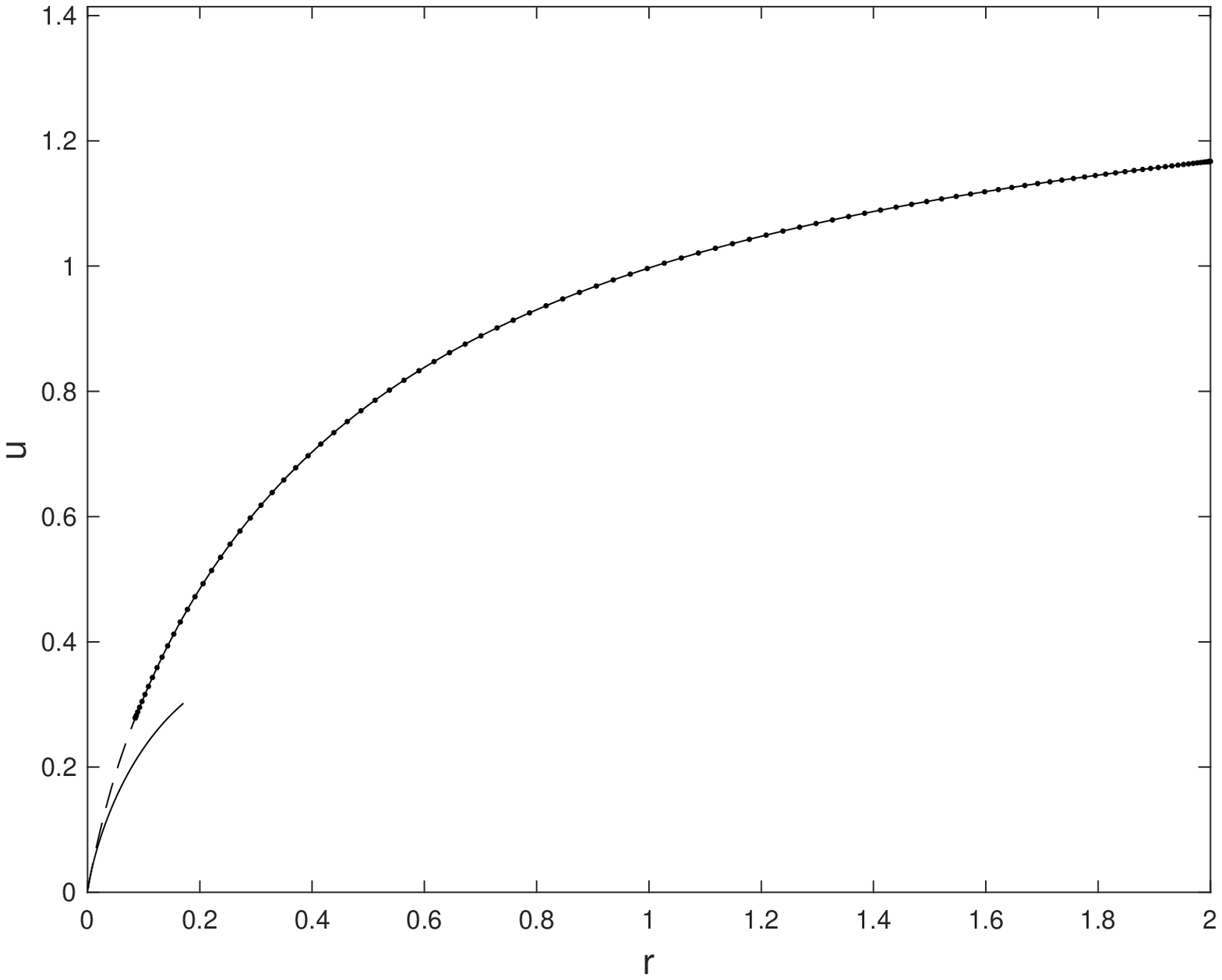}}
	\caption{The function $T(\sigma)$ with the asymptotic information and interpolating the computed data with the 10 leftmost asymptotic points included in the data set.}
	\label{fig:Tinterp}
\end{figure}
\begin{figure}[h!]
	\centering
	\scalebox{0.45}{\includegraphics{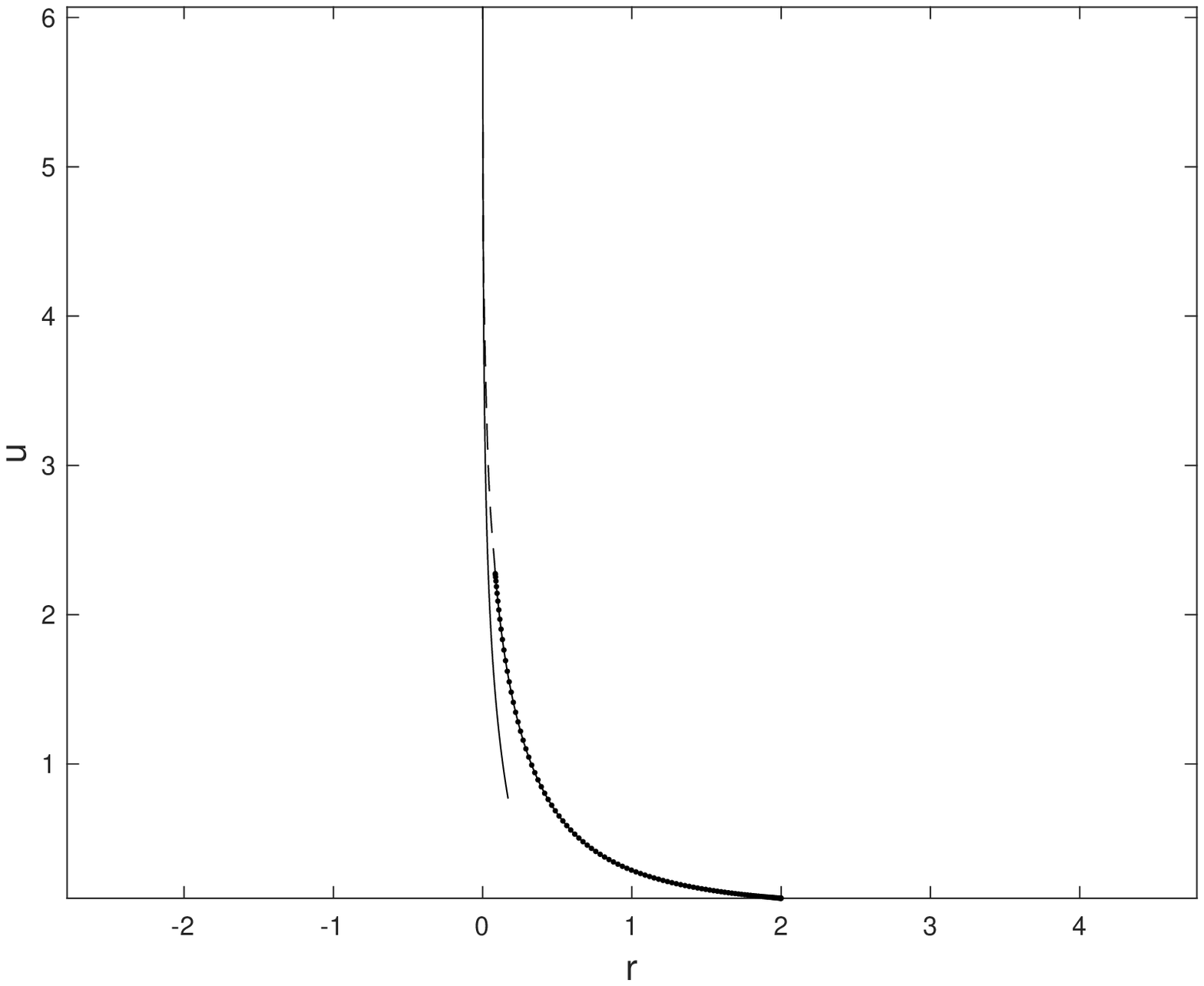}}
	\caption{The function $T^\prime(\sigma)$ with the asymptotic information and interpolating the computed data with the 10 leftmost asymptotic points included in the data set.}
	\label{fig:DTinterp}
\end{figure}

Finally, we make an effort to see better what is happening as $\sigma\rightarrow 0$.  In his Theorem~3.2, Turkington \cite{Turkington} showed that 
\begin{equation}\label{eqn:Turkest}
T(\sigma) \sim -\sigma \log(\sigma).
\end{equation}
Though it is not immediately clear when this asymptotic estimate begins to become a good approximation for the function $T(\sigma)$, we will use this to approximate the derivative:
\begin{equation}\label{eqn:DTurkest}
T^\prime(\sigma) \sim - \log(\sigma) - 1.
\end{equation}
In Figure~\ref{fig:Tinterp} we show the results of the algorithm described above and we include two new curves.  The (lower) dotted curve in that figure is the estimate \eqref{eqn:Turkest} plotted from the radius 0.00085 to the radius 0.170.  Clearly the asymptotic estimate is not very good near 0.1.  Given the theme of the numerical portions of this paper, we plotted the dotted curve based on 100 Chebyshev points over that interval.  Then we took the leftmost 10 of these points and the corresponding data for $T(\sigma)$ there and included this into an array with the data points we had from our algorithm to compute $T(\sigma)$.  Then we interpolated this data with a cubic spline at 100 Chebyshev points between 0.00085 and 0.085.  The result is graphed with a dashed curve in the figure.  This gives an estimate for the curve $T(\sigma)$ in this range, though we do not offer any comment as to how accurate it might be.

The more immediately important estimation process is to treat $T^\prime(\sigma)$ using \eqref{eqn:DTurkest} in the same fashion.  We follow the exact same steps just described for $T(\sigma)$, though here we use the data we computed for $T^\prime(\sigma)$ using our algorithm, and the formula \eqref{eqn:DTurkest} for the asymptotic data.  The result is graphed in Figure~\ref{fig:DTinterp}.

\begin{figure}[t]
	\centering
	\scalebox{0.45}{\includegraphics{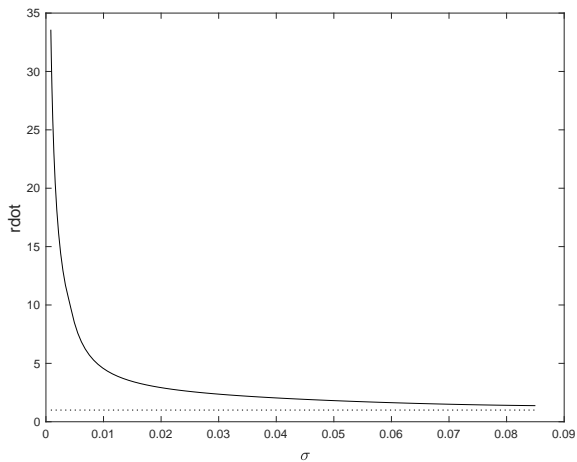}}
	\caption{The endpoint $\dot r(0,\sigma)$ is graphed as a function of $\sigma$.  The lowest part of the curve $\dot r(\phi,\sigma)$ is at $\phi = \pi/2$, and this is graphed with a dotted line. }
	\label{fig:minrdotinterp}
\end{figure}

We are now able to repeat our numerical experiment on this interpolated data.   We again construct numerical approximations of the solutions of \eqref{eqn:dotr}-\eqref{eqn:dotu} with the boundary conditions \eqref{eqn:bcdotr}-\eqref{eqn:bcdotu} where we need to append that system with the original \eqref{eqn:drdphi}-\eqref{eqn:dudphi} where $r(\pi/2,\sigma) = \sigma$ and $u(\pi/2,\sigma) = T(\sigma)$.  Then we use Matlab's ode45 to solve this system from $\phi = \pi/2$ to $\phi= 0$ to get the ``top'' portion of the solution.  Here we  ask for 11 digits of accuracy in both the absolute and relative error.  We then use all 100 values of $\sigma$ considered in this interpolation process to graph the endpoint $\dot r(0,\sigma)$  as a function of $\sigma$.  In all cases here, this endpoint is not the  lowest part of the curve $\dot r(\phi,\sigma)$, however  the lowest part of the curve occurs at $\phi = \pi/2$, and this is also included in Figure~\ref{fig:minrdotinterp}.

 These results strongly supports the result claimed in \cite{ENS} by numerically satisfying one of the conditions in Theorem~\ref{thm:hyp} by showing $\dot r(\phi_0) > 0$ for $\phi_0\in[0,\pi/2]$ in all computed cases.

\FloatBarrier

\end{document}